\documentclass[12pt,titlepage]{article}
\usepackage{amssymb}
\input epsf
\newcommand{\proof}{{\noindent \bf Proof. }}
\newtheorem{thm}{Theorem}
\newtheorem{defi}{Definition}

\newtheorem{lem}{Lemma}[section]
\newtheorem{cor}[thm]{Corollary}
\newtheorem{prop}[thm]{Proposition}

\newcommand{\HH}{{\cal H}}

\def\2{\mathbb Z_2}
\def\S{\mathbb S}
\def\ind{{\rm ind}}
\def\coind{{\rm coind}}
\def\susp{{\rm susp}}

\newcommand\mbf[1]{\mbox{\boldmath$#1$}}

\title{Local chromatic number and distinguishing the strength of topological
obstructions} 

\author{
      {\bf G\'abor Simonyi}\thanks{Research partially supported by the
Hungarian Foundation for Scientific Research Grant (OTKA) Nos.\ T037846,
T046376, AT048826, and NK62321.}\\
 Alfr\'ed R\'enyi Institute of Mathematics\\ 
 Hungarian Academy of Sciences\\
 1364 Budapest, POB 127, Hungary\\
 {\tt simonyi@renyi.hu}
\and
   {\bf G\'abor Tardos}\thanks{Research partially supported by the
NSERC grant 611470 and the
Hungarian Foundation for Scientific Research Grant (OTKA) Nos.\ T037846,
T046234, AT048826, and NK62321.}\\
 School of Computing Science, Simon Fraser University \\
 Burnaby BC, Canada V5A 1S6\\
 and\\
 Alfr\'ed R\'enyi Institute of Mathematics \\
 Hungarian Academy of Sciences\\
 1364 Budapest, POB 127, Hungary\\
 {\tt tardos@cs.sfu.ca}
\and
   {\bf Sini\v{s}a T. Vre\'cica}\thanks{Supported by the Serbian
Ministry of Science, Grant 144026.} \\ 
Faculty of Mathematics, University of Belgrade \\
Studentski trg 16, P.O.B. 550 \\
11000 Belgrade, Serbia\\ 
{\tt vrecica@matf.bg.ac.yu}} 

\date{}

\begin{document}
\maketitle

\begin{abstract}
The local chromatic number of a graph $G$ is the number of colors appearing in
the most colorful closed neighborhood of a vertex minimized over all proper
colorings of $G$. We show that two specific topological obstructions that have
the 
same implications for the chromatic number have different implications for the
local chromatic number. These two obstructions can be formulated in terms of
the homomorphism complex ${\rm Hom}(K_2,G)$ and its suspension, respectively. 

These investigations follow the line of research initiated by Matou\v{s}ek and
Ziegler who recognized a hierarchy of the different topological expressions
that can serve as lower bounds for the chromatic number of a graph.

Our results imply that the local chromatic number of $4$-chromatic Kneser,
Schrijver, Borsuk, and generalized Mycielski graphs is $4$, and more
generally, that $2r$-chromatic versions of these graphs have local chromatic
number at least $r+2$. This lower bound is tight in several cases by results
in \cite{ST}. 
\bigskip
\bigskip
\bigskip
\par
\noindent
{\it 2000 Mathematics Subject Classification.} Primary 05C15; Secondary 57M15. 
\par
\noindent
{\it Key words and phrases}. Local chromatic number, box complex, Borsuk-Ulam
theorem.
\end{abstract}

\section{Introduction}

The local chromatic number is a coloring type graph parameter defined by
Erd\H{o}s, F\"uredi, Hajnal, Komj\'ath, R\"odl, and Seress \cite{EFHKRS} in
1986. It is the number of colors appearing in the most colorful
closed neighborhood of a vertex minimized over all proper colorings of the
graph. Using the notation $N(v)=N_G(v):=\{w: vw\in E(G)\}$, the formal
definition is as follows.  

\begin{defi} \label{defi:lochr} {\rm(\cite{EFHKRS})}
The {\em local chromatic number} $\psi(G)$ of a graph $G$ is
$$\psi(G):=\min_c \max_{v\in V(G)} |\{c(u): u \in N(v)\}|+1,$$
where the minimum is taken over all proper colorings $c$ of $G$.
\end{defi}

Considering closed neighborhoods $N(v)\cup\{v\}$ results in a simpler form
of the relations with other coloring parameters and explains the +1 term in
the definition. 

It is clear that $\psi(G)$ is always bounded from above by $\chi(G)$, the
chromatic number of $G$. It is also easy to see that $\psi(G)=2$ is equivalent
to $\chi(G)=2$. However, as it is proven in \cite{EFHKRS}, cf.\ also
\cite{Fur}, there exist graphs with $\psi(G)=3$ and $\chi(G)$ arbitrarily
large. In this sense the local chromatic number is highly independent of the
chromatic number. 

On the other hand, it was observed in \cite{KPS} that the fractional chromatic
number $\chi_f(G)$ serves as a lower bound, i.e., $\chi_f(G)\leq \psi(G)$
holds. (For the definition and basic properties of the fractional chromatic
number we refer to the books \cite{SchU} and \cite{GR}.) This motivated in
\cite{ST} the study of the local chromatic number of graphs that have a large
gap 
between their ordinary and fractional chromatic numbers. Basic examples of such
graphs include Kneser graphs and Mycielski graphs (see \cite{SchU}) and their
variants, the so-called Schrijver graphs (see \cite{Mat}, \cite{Schr}) and
generalized Mycielski graphs (see \cite{GyJS}, \cite{Mat}, \cite{Stieb},
\cite{Tar}).  
Another common feature of these graphs is that their chromatic number is (or
at least can be) determined by the topological method initiated by Lov\'asz in
\cite{LLKn}. In \cite{ST} it is proved that for all these graphs of chromatic
number $t$ one has
$$\psi(G)\ge\left\lceil t\over 2\right\rceil+1,$$ 
and showed several cases when this bound is tight. In all those cases,
however, we have an odd $t$, in particular, the smallest chromatic number for
which we have shown some Schrijver graphs, say, with smaller local than
ordinary chromatic number is $5$, in spite of the fact, that the lower bound
$\lceil{t\over 2}\rceil+1$ is smaller than $t$ already for $t=4$. 
In this paper we show that whether $t=4$ or $5$ is optimal in the above sense
depends on the particular topological method that gives the chromatic number
of the graph. An analogous difference between the best possible lower bound on
the local chromatic number will be shown to exist for $2r$-chromatic
graphs in general. 
In \cite{ST} two possible topological requirements were considered that 
make the chromatic number of a graph at least $t$. Here we show in one hand
that since the graphs mentioned above satisfy the stronger of these two
requirements, they also satisfy $\psi\ge r+2$ in the $t=2r$ case. On the
other hand, we show that the general lower bound in \cite{ST},
which is derived from the weaker topological requirement considered,
is tight in the sense that for all $t$ there exist graphs for which the
above lower bound applies with equality. In particular, this
shows that the two kinds of topological obstructions for graph coloring have
different implications in terms of the local chromatic number. This
consequence is in the spirit of the investigations by 
Matou\v{s}ek and Ziegler \cite{MZ} about the hierarchy they discovered among
the different topological techniques bounding the chromatic number. 

Some of the results (concerning the case $t=4$) below were announced in the
brief summary \cite{STB}. 

\section{Preliminaries}

\subsection{Topological preliminaries} \label{topre}

The following is a brief overview of some of the topological concepts we
need. We refer to \cite{Bjhand,Hat} and \cite{Mat} for basic
concepts and also for a more detailed discussion of the notions and facts
given below. 
We use the notations of \cite{Mat}.
\smallskip

A {\em $\mathbb{Z}_2$-space} (or {\em involution space}) is a pair $(T,\nu)$
of a topological space $T$ and the involution $\nu:T\to T$, which is
continuous and satisfies that $\nu^2$ is the identity map. The points $x\in T$
and $\nu(x)$ are called {\em antipodal}. The involution
$\nu$ and the $\mathbb{Z}_2$-space $(T,\nu)$ are {\em free} if
$\nu(x)\ne x$ for all points $x$ of $T$. If the involution is understood
from the context we speak about $T$ rather than the pair $(T,\nu)$. This is
the case, in particular, for the unit sphere $\mathbb S^d$ in ${\mathbb
  R}^{d+1}$ 
with the 
involution given by the central reflection ${\mbf x}\mapsto-{\mbf x}$. 
A continuous map $f:S\to T$ between $\mathbb{Z}_2$-spaces $(S,\nu)$ and
$(T,\pi)$ is a {\em $\mathbb{Z}_2$-map} (or an {\em equivariant map}) if it
respects the respective
involutions, that is $f\circ\nu=\pi\circ f$. If such a map exists we write
$(S,\nu)\to(T,\pi)$. If
$(S,\nu)\to(T,\pi)$ does not hold we write $(S,\nu)\not\to(T,\pi)$. If both
$S\to T$ and $T\to S$ we call the $\mathbb Z_2$-spaces $S$ and $T$ \ $\mathbb
Z_2$-equivalent and write $S\leftrightarrow T$.

We sometimes refer to homotopy equivalence and $\2$-homotopy equivalence 
(i.e., homotopy equivalence given by $\2$-maps), but will use only 
the following two simple 
observations. First, if the $\2$-spaces $S$ and $T$ are $\2$-homotopy
equivalent, then $S\leftrightarrow T$. Second, if the space $S$ is
homotopy equivalent to a sphere $\mathbb S^h$ (this relation is between
topological 
spaces, not $\2$-spaces), then $S$ is $(h-1)$-connected and therefore ${\mathbb
S}^h\to(S,\nu)$ for any involution $\nu$, cf. \cite{Mat} (proof of
Proposition 5.3.2 (iv), p. 97). In the other direction we have
$(S,\nu)\rightarrow \mathbb S^h$ if $(S,\nu)$ is the body of a $h$-dimensional
free simplicial $\2$-complex. (See below the definition of the latter.)

\medskip

The $\mathbb{Z}_2$-index of a $\mathbb{Z}_2$-space $(T,\nu)$ is defined 
(see e.g.\ \cite{MZ,Mat}) as 
$$\ind(T,\nu):=\min\{d\ge 0:(T,\nu)\to \mathbb S^d\},$$
where $\ind(T,\nu)$ is set to be $\infty$ if
$(T,\nu)\not\to \mathbb S^d$ for all $d$. 
\smallskip
\par
\noindent
The $\mathbb{Z}_2$-coindex of a $\mathbb{Z}_2$-space $(T,\nu)$ is defined as
$$\coind(T,\nu):=\max\{d\ge 0:\mathbb S^d\to(T,\nu)\}.$$
If such a map exists for all $d$, then we set $\coind(T,\nu)=\infty$. 
Thus, if $(T,\nu)$ is not free, we have ${\rm
ind}(T,\nu)=\coind(T,\nu)=\infty$. 
\medskip
\par
\noindent
Note that $S\to T$ implies $\ind(S)\le\ind(T)$ and $\coind(S)\le\coind(T)$. In
particular, $\2$-equivalent spaces have equal index and also equal coindex.

The celebrated Borsuk-Ulam Theorem can be stated in many equivalent
forms. Here we state four of them. For more equivalent versions and several
proofs we refer to \cite{Mat}. Here (i)--(iii) are all standard forms of the
Borsuk-Ulam Theorem, while (iv) is clearly equivalent to (iii).
\medskip

\noindent
{\bf Borsuk-Ulam Theorem.}
{\em\begin{description}
\item[(i)] For every continuous map $f: \mathbb S^k\to\mathbb R^k$ there 
exists ${\mbf x}\in \mathbb S^k$ for which $f({\mbf x})=f(-{\mbf x})$.
\item[(ii)] (Lyusternik-Schnirel'man version) Let $d\ge0$ and let $\HH$ be a
collection 
of open (or closed) sets covering $\mathbb S^d$ with no $H\in\HH$ containing a pair of
antipodal points. Then $|\HH|\ge d+2$.
\item[(iii)] $\mathbb S^{d+1}\not\to \mathbb S^d$ for any $d\ge 0$.
\item[(iv)] For a $\mathbb{Z}_2$-space $T$ we have $\ind(T)\ge
\coind(T)$.  
\end{description}}
\medskip
\par
\noindent

The suspension $\susp(S)$ of a topological space $S$ is defined as the factor
of the space $S\times[-1,1]$ that identifies all the points in $S\times\{-1\}$
and identifies also the points in $S\times\{1\}$. If $S$ is a $\2$-space with
the involution $\nu$, then the suspension $\susp(S)$ is also a $\2$-space with
the involution $(x,t)\mapsto(\nu(x),-t)$. Any $\2$-map $f:S\to T$ naturally
extends to a $\2$-map $\susp(f):\susp(S)\to\susp(T)$ given by
$(x,t)\mapsto(f(x),t)$. We have $\susp(\mathbb S^n)\cong \mathbb S^{n+1}$ with a
$\2$-homeomorphism. These observations show the well known inequalities below.

\begin{lem}\label{suspension}
For any $\2$-space $S$ \ $\ind(\susp(S))\le\ind(S)+1$ and
$\coind(\susp(S))\ge\coind(S)+1$.
\end{lem}

A(n abstract) {\em simplicial complex} $K$ is a non-empty, hereditary set
system. In this paper we consider only finite simplicial complexes.
The non-empty sets in $K$ are called {\em simplices}. The {\em dimension} of a
$\sigma\in K$ is $\dim(\sigma)=|\sigma|-1$. A simplex of dimension $k$
is called a {\em $k$-simplex}. The {\em dimension} of $K$ is defined as
$\max\{\dim(\sigma):\sigma\in K\}$. We call the set
$V(K)=\{x:\{x\}\in K\}$
the set of {\em vertices} of $K$. In a {\em geometric realization} of $K$ a
vertex $x$ corresponds to a point $||x||$ in a Euclidean space, a simplex
$\sigma$ corresponds to its {\em body}, the convex hull of its vertices:
$||\sigma||={\rm conv}(\{||x||:x\in\sigma\})$. We assume that the points
$||x||$ for $x\in\sigma$ are affine independent, and so $||\sigma||$ is a
geometric simplex. We also assume that disjoint simplices have disjoint
bodies. The body of the complex $K$ is $||K||=\cup_{\sigma\in K}||\sigma||$.
$||K||$ is determined up to homeomorphism by $K$. Any point in $p\in||K||$
has 
a unique representation as a convex combination $p=\sum_{x\in
V(K)}\alpha_x||x||$ such that $\{x:\alpha_x>0\}\in K$.

A {\em simplicial map} $f:K\to L$ maps the vertices of a simplicial complex
$K$ to the vertices of another simplicial complex $L$ such that the image of
a simplex of $K$ is a simplex in $L$. Such a map can be linearly
extended to the bodies of all simplices in $K$ giving a continuous map
$||f||:||K||\to||L||$. A simplicial complex with a simplicial
involution is called a {\em simplicial $\2$-complex}. 

The barycentric subdivision $sd(K)$ of a simplicial complex $K$ is the family
of chains (subsets linearly ordered by inclusion) of simplices of $K$. The
standard geometric realization (each simplex is represented by a point in its
relative interior) gives $||sd(K)||=||K||$.

\subsection{Topological lower bounds on the chromatic number}

The topological method for bounding the chromatic number can be described by
the following scheme. One assigns a $\2$-space to all graphs in such a
way that whenever a homomorphism from $F$ to $G$ exists this implies the
existence of a $\2$-map from the space assigned to $F$ to that assigned to
$G$. Colorability with $m$ colors is equivalent to the existence of a
homomorphism to $K_m$. If one shows that no $\2$-map exists from the space
assigned to $G$ to the space assigned to $K_m$, then it proves that $G$ is not
$m$-colorable. In the cases we consider the space assigned to $K_m$ will be
$\2$-homeomorphic to 
$\mathbb S^{f(m)}$ with $f(m)=m-2$ or $m-1$ depending on which of the two
space assignments discussed below is used. 
Thus if $G$ is $m$-colorable, then the $\2$-index of the space assigned
to $G$ must not be more than $f(m)$. If it is more than $f(m)$ that implies
$\chi(G)>m$. Thus
we can bound the chromatic number from below by giving a lower 
bound on the index of a certain $\2$-space. This is often done by
actually bounding its coindex from below. By the Borsuk-Ulam theorem (form
(iv)) this also provides a lower bound on the index. 

One way to assign a $\2$-space to a graph $G$ is via defining some
simplicial complex, a so-called box complex, and
considering the body of this complex. Following the papers \cite{AFL,Kriz}
Matou\v{s}ek and Ziegler \cite{MZ} defines several box complexes that turn out
to fall into two categories in the sense that their index (or coindex) assumes
one of only two values. (This is proven in \cite{MZ}, but Csorba \cite{Cs} and
\v{Z}ivaljevi\'c \cite{Ziv} gives further explanation of this fact by showing
that the homotopy type of all these complexes is one of only two different
kinds.) 
One representative of both of these types are given in the two definitions
below. (In the second
case, for simplicity, we speak about a cell complex and its body as the
corresponding topological space. It
is also $\2$-homotopy equivalent to some of the known box complexes as
remarked after Definition~\ref{hom}.) 

For subsets
$S,T\subseteq V(G)$ we denote the set $S\times\{1\}\cup
T\times\{2\}$ by $S\uplus T$.
For $v\in V(G)$ we denote by $+v$ the vertex
$(v,1)\in\{v\}\uplus\emptyset$ and $-v$ denotes the vertex
$(v,2)\in\emptyset\uplus\{v\}$.

\begin{defi}
The {\em box complex} $B_0(G)$ is a simplicial complex on the vertices
$V(G)\times\{1,2\}$. For subsets $S,T\subseteq V(G)$ the set $S\uplus
T:=S\times\{1\}\cup T\times\{2\}$ forms a simplex if and only if 
$S\cap T=\emptyset$ and the complete bipartite graph with sides $S$ and $T$ is
a subgraph of $G$. The simplicial involution switching $+v$ and $-v$ for $v\in
V(G)$ makes $B_0(G)$ a simplicial $\2$-complex and $||B_0(G)||$ a free
$\mathbb{Z}_2$-space. 
\end{defi}

Note that $V(G)\uplus\emptyset$ and $\emptyset\uplus V(G)$ are simplices of
$B_0(G)$. 

\begin{defi} \label{hom}
The {\em hom space} $H(G)$ of $G$ is the subspace of
$||B_0(G)||$ consisting of
those points $p\in||B_0(G)||$ that, when written as a convex combination
$p=\sum_{x\in V(B_0(G))}\alpha_x||x||$ with $\{x:\alpha_x>0\}\in B_0(G)$ give
$\sum_{x\in V(G)\uplus\emptyset}\alpha_x=1/2$. This space can also be
considered as the body of a cell complex as follows. Let the {\em hom complex}
${\rm Hom}(K_2,G)$ of $G$ be the cell complex with cells $S\uplus T\in B_0(G)$
with 
$S\ne\emptyset\ne T$. We call $S\uplus T\in{\rm Hom}(K_2,G)$ a {\em cell} of
the 
complex and $||S\uplus T||\cap H(G)$ is the {\em body} of this cell. The
vertices of ${\rm Hom}(K_2,G)$ are of the form $\{x\}\uplus\{y\}$ with
$\{x,y\}\in E(G)$.

We consider ${\rm Hom}(K_2,G)$ as a $\2$-complex and $H(G)$ as a $\2$-space
with the involution inherited from $B_0(G)$.  
\end{defi}

The cell complex ${\rm Hom}(K_2,G)$ is a special case of the more general
homomorphism complexes ${\rm Hom}(F,G)$, 
see \cite{BK}. 
The hom space $H(G)$ can also be considered as the body of a
simplicial complex $B_{chain}(G)$, where $B_{chain}(G)$ is the first
barycentric subdivision of ${\rm Hom}(K_2,G)$, see \cite{MZ}. 
The latter is also $\2$-homotopy equivalent to another simplicial box
complex $B(G)$ (for a formal definition of $B(G)$, cf.\ \cite{MZ}) where
$B(G)$ is the hereditary closure of ${\rm Hom}(K_2,G)$ and it
differs from $B_0(G)$ only by not containing those simplices $S\uplus T$ where
the elements of one of the sets $S$ and $T$ do
not have a common neighbor in $G$ (implying emptiness of the other set).

A useful connection between $B_0(G)$ and $H(G)$ follows from 
results of Csorba. Namely, Csorba \cite{Cs} proves the $\2$-homotopy
equivalence of $||B_0(G)||$ and the suspension of the body of the other box
complex $B(G)$ mentioned above. Further, he proves, cf.\ also
\v{Z}ivaljevi\'c \cite{Ziv}, the $\2$-homotopy equivalence of $||B(G)||$ and
$H(G)$. (A weaker version of the latter equivalence, which already implies
the proposition below also follows from the results in \cite{MZ}.)  

\begin{prop}\label{csorba} {\rm (\cite{Cs,MZ,Ziv})}
$||B_0(G)||\leftrightarrow\susp(H(G))$.
\end{prop}

The box complex $B_0(K_m)$ is the
boundary complex of the $m$-dimensional {\em cross-polytope} (i.e., the
convex hull of the basis vectors and their negatives in ${\mathbb R}^m$), thus
$||B_0(K_m)||\cong \mathbb S^{m-1}$ with a $\2$-homeomorphism and
$\coind(||B_0(G)||)\le\ind(||B_0(G)||)\leq m-1$ is necessary for $G$ being
$m$-colorable. Similarly, $\coind(H(G))\le\ind(H(G))\leq m-2$ is also
necessary for $\chi(G)\leq m$ since $H(K_m)$ can be obtained
from intersecting the boundary of the $m$-dimensional cross-polytope with the
hyperplane $\sum x_i=0$, and therefore $H(K_m)\cong \mathbb S^{m-2}$ with a
$\2$-homeomorphism.  These four lower bounds on $\chi(G)$ can be arranged in a
single line of inequalities using Lemma~\ref{suspension} and
Proposition~\ref{csorba}:  
\begin{equation}\label{eq:chib1}
\chi(G)\ge\ind(H(G))+2\ge\ind(||B_0(G)||)+1\ge\coind(||B_0(G)||)+1\ge\coind(H(G))+2 
\end{equation}

The first two of the lower bounds to $\chi(G)$ above are (equivalent to) the
two strongest lower bounds in Matou\v{s}ek and Ziegler's Hierarchy Theorem
\cite{MZ}. We are able to say more on the last two bounds that were singled out
by the following definition in \cite{ST}. 

\begin{defi} \label{defi:topres}
We say that a graph $G$ is {\em topologically $t$-chromatic} if
$$\coind(||B_0(G)||)\ge t-1.$$
We say that a graph $G$ is {\em strongly topologically $t$-chromatic} if
$$\coind(H(G))\ge t-2.$$
\end{defi}

Note that if a graph is strongly topologically $t$-chromatic, then it is also 
topologically $t$-chromatic, and if $G$ is topologically $t$-chromatic, then
$\chi(G)\ge t$.

Examples of strongly topologically $t$-chromatic graphs are provided by
$t$-chromatic Kneser graphs, Schrijver graphs, generalized Mycielski
graphs. (For the formal definition of all these graphs,
see, e.g., \cite{Mat}, or \cite{ST}.)
One way to show that these graphs are strongly topologically $t$-chromatic is
to refer to another simplicial complex, the neighborhood complex ${\cal N}(G)$
of the graph $G$, introduced by Lov\'asz in \cite{LLKn}. Proposition 4.2 in
\cite{BK} states that $||{\cal N}(G)||$ is homotopy equivalent to $H(G)$ for
every graph $G$ (note that $||{\cal N}(G)||$ is {\em not} a $\2$-space, thus
this cannot be a $\2$-homotopy equivalence). Thus if ${\cal
N}(G)$ is homotopy equivalent to the sphere $\mathbb S^{t-2}$  then, by the
above 
result in \cite{BK} and the corresponding remark in the introductory part of
Subsection~\ref{topre}, we have $\coind(H(G))\ge t-2$. 
(In fact, since $H(G)$ is free, we have equality here.) For $t$-chromatic
Schrijver graphs Bj\"orner and de Longueville \cite{BjLo} proved that their
neighborhood complex is homotopy equivalent 
to $\mathbb S^{t-2}$. As Schrijver graphs are induced subgraphs of Kneser
graphs with 
the same chromatic number this proves strong topological $t$-chromaticity for
both $t$-chromatic Kneser graphs and Schrijver graphs. An analogous result
about the homotopy equivalence of the neighborhood complex of $t$-chromatic
generalized Mycielski graphs and $\mathbb S^{t-2}$ was proved by Stiebitz
\cite{Stieb}, cf.\ also \cite{GyJS} and \cite{Mat}. 
There is a similar result due to Lov\'asz \cite{LLgomb} for a finite subgraph
of the Borsuk graph $B(t-1,\alpha)$ (see Definition~\ref{defi:Bogr}) that 
we will return to in the proof of Lemma~\ref{lem:equiB}.
We remark that the strong topological $t$-chromaticity of $t$-chromatic Kneser
graphs and Schrijver graphs can also be seen more directly from the results of
B\'ar\'any \cite{Bar} and Schrijver \cite{Schr}. For more details about this, 
cf.\ Proposition 8 in \cite{ST}.  

For examples of graphs that are topologically $t$-chromatic
but not strongly topologically $t$-chromatic we refer to the detailed
discussion in Sections~\ref{versus} and \ref{remark14}. A longer list of
topologically $t$-chromatic graphs is given in \cite{ST3}.

\section{Local chromatic number and covering the sphere}
\label{bounds}

In \cite{ST} the following lower bound on the local chromatic number
of topologically $t$-chromatic graphs is proved. 

\begin{thm} \label{thm:lowb} {\rm (\cite{ST})} 
If $G$ is topologically  $t$-chromatic for some $t\ge2$, then 
$$\psi(G)\ge\left\lceil t\over 2\right\rceil+1.$$
\end{thm}

The proof was based on an old topological theorem of Ky Fan \cite{kyfan} which
generalizes the Borsuk-Ulam theorem. 
It was also shown in \cite{ST} that this lower bound is tight for several
Schrijver graphs, generalized Mycielski graphs, and Borsuk graphs of odd
chromatic number. 

Here we prove a similar but somewhat different lower bound than the one in 
Theorem~\ref{thm:lowb}. It applies only for strongly topologically
$t$-chromatic graphs and gives the same conclusion if $t$ is odd, thus it is a
weaker statement in that case. For $t$ even, however, the conclusion is also
slightly stronger. 

\begin{thm} \label{thm:ps4}
If a graph $G$ is strongly topologically $t$-chromatic for $t\ge 3$,
then $$\psi(G)\ge \left\lfloor{t\over 2}\right\rfloor +2.$$ 
\end{thm}

To prove that a similar statement is not true for topologically $t$-chromatic
graphs we will show in Section~\ref{versus} for every $r\ge 2$ a
topologically $2r$-chromatic graph $G$ with $\psi(G)=r+1$. 
By Theorem~\ref{thm:ps4} this
graph cannot be strongly topologically $2r$-chromatic. Thus together with
Theorem~\ref{thm:ps4} it 
proves that topological $t$-chromaticity and strong topological
$t$-chromaticity have different implications for the local chromatic
number. 

First we translate the problem into one concerning open covers of the sphere.

\begin{defi}\label{q} 
For a nonnegative integer parameter $h$ let $Q(h)$ denote the minimum number
$l$ for which $\mathbb S^h$ can be covered by open sets in such a way that no
point of 
the sphere is contained in more than $l$ of these sets and none of the covering
sets contains an antipodal pair of points.
\end{defi}

In the earlier paper \cite{ST} the first two authors arrived to the problem of
determining $Q(h)$ 
through local colorings of graphs. The same question was independently asked
by Micha Perles motivated by a related question of Matatyahu Rubin\footnote{We
are indebted to Imre B\'ar\'any \cite{b} and Gil Kalai \cite{GK} for this
information.}. After the
publication of \cite{ST} we learnt that this question was already
considered and settled in papers by \v{S}\v{c}epin \cite{Scep}, Izydorek,
Jaworowski \cite{IJ}, and
Jaworowski \cite{Jaworpre,Jaworgen}, cf. also Aarts and Fokkink \cite{AF}. The
$h=2$ case was solved even earlier by Shkliarsky \cite{Sklj}. (The papers
\cite{Scep}, \cite{IJ}--\cite{Jaworgen}
use the different but equivalent formulation that we will see in
Lemma~\ref{lem:equiB}~(v) below. This equivalence is already implicit in
\cite{Scep}, cf. also \cite{AF}.) 

\begin{thm}\label{qfalse} {\rm (\cite{Scep}, \cite{IJ}--\cite{Jaworgen})}
For every $h\ge 1$ 
$$Q(h)=\left\lfloor h\over
2\right\rfloor+2.$$
\end{thm}

\medskip
\par
\noindent
{\it Remark 1.}
The results in \cite{ST} had the implications $\left\lceil h\over
  2\right\rceil+1\le Q(h)\le\left\lfloor h\over 2\right\rfloor+2$, where the
lower bound followed from Ky Fan's theorem \cite{kyfan}, cf. Corollary 18 in
\cite{ST}. The slightly stronger Corollary 17 of \cite{ST} can also be read out
from results in \cite{AF}.   
\hfill$\Diamond$

The relevance of the value of $Q(h)$ to local colorings will be
clarified in Lemma \ref{lem:equiB} below. 
One of the conditions in the lemma uses the concept of Borsuk graphs. 
Their appearance in the equivalent conditions for $Q(h)\leq l$ parallels the
fact that the Borsuk-Ulam theorem is equivalent to stating the chromatic
number of Borsuk graphs (of appropriate parameters) as remarked by Lov\'asz in
\cite{LLgomb}.

\begin{defi} \label{defi:Bogr}
The Borsuk graph $B(n,\alpha)$ of parameters $n$ and $0<\alpha<2$ is the
infinite graph whose vertices are the points of the unit sphere in ${\mathbb
R}^n$ (i.e., $\mathbb S^{n-1}$) and whose edges connect the pairs of points
with distance at least $\alpha$. 
\end{defi}

\begin{lem} \label{lem:equiB}
The following five statements are equivalent for every $h$ and $l$.
\begin{description}
\item[(i)] $Q(h)\leq l$, i.e., $\mathbb S^h$ can be covered by
open sets such that none of them contains an antipodal pair of points and no
${\mbf x}\in \mathbb S^h$ is contained in more than $l$ of these sets.
\item[(ii)] $\mathbb S^h$ can be covered by a finite number of closed sets
  such that 
none of them contains an antipodal pair of points and no ${\mbf x}\in \mathbb
S^h$ is 
contained in more than $l$ of these sets.
\item[(iii)] There exists $0<\alpha<2$ for which $\psi(B(h+1,\alpha))\le l+1$. 
\item[(iv)] There exists a finite graph $G$ with $\coind(H(G))\ge h$ (i.e.,
a strongly topologically $(h+2)$-chromatic graph) such that $\psi(G)\le l+1$.
\item[(v)] There is a continuous map $g$ from $\mathbb S^h$ to the body
  $||K||$ of a 
finite simplicial complex $K$ of dimension at most $l-1$ satisfying $g({\mbf
x})\neq g(-{\mbf x})$ for all ${\mbf x}\in \mathbb S^h$.
\end{description}
\end{lem}

We note that, as already mentioned, the equivalence of (ii) and (v) is already
implicit in \cite{Scep} and is also contained partially in Lemma 5 of
\cite{AF}.  

We also note that for a finite graph $G$ the property 
$\coind(H(G))\ge h$ can also be described in terms of Borsuk graphs: it is
equivalent to the existence of a homomorphism from $B(h+1,\alpha)$ to $G$ for
appropriately large $\alpha<2$, cf. \cite{ST}. 
\smallskip
\par
\noindent
\proof

\noindent
(ii)$\Rightarrow$(iii):
Consider a covering $\cal A$ as in (ii). Consider the closed sets in the
covering as colors and color each point of
$\mathbb S^h$ with one of the sets containing it. We need to prove that if
$\alpha<2$ 
is large enough this is a proper coloring establishing $\psi(B(h+1,\alpha))\le
l+1$.

We may assume that $|{\cal A}|>l$, otherwise we can add singleton sets. For
each ${\mbf x}\in \mathbb S^h$ let $g({\mbf x})$ be the $(l+1)^{st}$ smallest
distance 
of a set $A\in\cal A$ from ${\mbf x}$. Since $g$ is the $(l+1)^{st}$ level
of a finite set of continuous functions, $g$ is continuous. Since $\mathbb
S^h$ is 
compact, $g$ attains its minimum $g({\mbf x}_0)$. Since the covering sets are
closed and ${\mbf x}_0$ is contained in at most $l$ of them,  
$g({\mbf x}_0)>0$. For any set $A\in\cal A$ the disjoint sets $A$
and $-A$ are compact and thus they have a positive distance. Let $\delta>0$ be
smaller than the minimum of $g$ and also smaller than the distance between
$A$ and $-A$ for all the sets $A\in\cal A$. We choose
$\alpha=\sqrt{4-\delta^2}$. With this choice the vertex $\mbf x$ of
$B(h+1,\alpha)$ is connected to the vertex $\mbf y$ exactly if the distance
between $\mbf y$ and $-\mbf x$ is at most $\delta$.

Let $\mbf x$ be a vertex of the Borsuk graph of color $A\in\cal A$. Any vertex
$\mbf y$ connected to $\mbf x$ is closer to $-\mbf x$ and hence to $-A$ then
$\delta$, therefore it cannot be contained in $A$. This shows that the
coloring is proper.

Consider the colors of the neighbors of $\mbf x$. These are sets with distance
at most $\delta$ from $-\mbf x$. From $g(-{\mbf x})>\delta$ it follows that the
number of these colors is at most $l$ as claimed.
\smallskip

\noindent
(iii)$\Rightarrow$(iv): Lov\'asz gives in \cite{LLgomb} a finite graph
$G_P\subseteq B(h+1,\alpha)$ which has the
property that its neighborhood complex ${\cal N(G)}$ is homotopy equivalent to
$\mathbb S^h$. Proposition 4.2 in \cite{BK} states that ${\cal N}(F)$ is
homotopy 
equivalent to $H(F)$ for every graph $F$, thus $\coind(H(G_P))\ge h$. As
$G_P\subseteq B(h+1,\alpha)$ we have $\psi(G_P)\le\psi(B(h+1,\alpha))\le l+1$.
\smallskip

\noindent
(iv)$\Rightarrow$(i): 
Consider a proper coloring $c$ of $G$ achieving $\psi(G)\leq l+1$ and let $m$
be the  number of colors used. First we give an at most $l$-fold covering of
$H(G)$ by open sets $U_1,\dots, U_m$. Let $y\in H(G)$ and let
$Z_y\uplus T_y$ be the minimal cell of ${\rm Hom}(K_2,G)$ (or equivalently, the
minimal simplex of $B_0(G)$) whose body contains $y$. We
let $y$ belong to $U_i$ if and only if there is some vertex $v\in Z_y$ for
which $c(v)=i$. It  is clear that the sets $U_i$ obtained this way are
open. As $Z_y\ne\emptyset$ the point $y$ is covered by some $U_i$. As
$T_y$ is not empty, we can choose a vertex $w\in T_y$. All vertices $v\in Z_y$
are neighbors of $w$, so by the definition of $\psi(G)$ these vertices have at
most $l$ different colors. Therefore $y$ is covered by at most $l$ sets $U_i$.
The sets $U_i$ therefore form an at most $l$-fold covering of $H(G)$. For
antipodal points $y,y'\in H(G)$ we have $Z_{y'}=T_y$. If
$y$ and $y'$ are contained in the same set $U_i$, then we find vertices $v\in
Z_y$ and $w\in T_y$ of the same color $i$. As $v$ and $w$ are adjacent and
$c$ is a proper coloring this is impossible, so the sets $U_i$ contain no
antipodal pairs of points.

By the condition $\coind(H(G))\ge h$ there is a $\mathbb{Z}_2$-map
$f:\mathbb S^h\to H(G)$. Now we define $A_i:=\{{\mbf x}\in \mathbb S^h:
f({\mbf x})\in 
U_i\}$. It is straightforward, that the open sets $A_1,\dots, A_m$ provide a
covering required. 

\smallskip

\noindent
(i) $\Rightarrow$ (v):
Assume that (i) holds. As $\mathbb S^h$ is compact we can assume that the open
cover 
is finite, it consists of the sets $A_1,\ldots,A_m$. Let $K$ be the simplicial
complex having vertices $[m]=\{1,\dots,m\}$ and all $l$-subsets of $[m]$ as
maximal simplices.
Define $g: \mathbb S^h\to||K||$ as follows. Let $d_i({\mbf x})$ be the
distance of 
${\mbf x}\in \mathbb S^h$ 
from $\mathbb S^h\setminus A_i$. Note that $d_i({\mbf
x})>0$ if and only if ${\mbf x}\in A_i$. We normalize $d_i$ to get
$\alpha_i({\mbf x})=d_i({\mbf x})/(\sum_{i=1}^m d_i({\mbf x}))$.  
Now set $g({\mbf x})$ to be the formal convex combination of the
vertices of  $K$ given by $\sum_{i=1}^m \alpha_i({\mbf x})||i||$. Since no
${\mbf x}\in \mathbb S^h$ is covered by 
more than $l$ of the sets $A_i$ the images are indeed in $||K||$. As the sets
$A_i$ do not contain antipodal points we have $g({\mbf x})\ne g(-{\mbf x})$,
furthermore the minimal simplices containing $g({\mbf x})$ and $g(-{\mbf x})$
are disjoint.
\smallskip

\noindent
(v) $\Rightarrow$ (ii): 
Let $g$ be a map as in (v). We assume that the minimal simplices containing
$g({\mbf x})$ and $g(-{\mbf x})$ are disjoint for every point ${\mbf x}\in
\mathbb S^h$. If this condition is
violated we consider an arbitrary geometric realization of $K$ and the
continuous function ${\mbf x}\mapsto{\rm dist}(g({\mbf x}),g(-{\mbf
x}))>0$. As $\mathbb S^h$ is compact this continuous function 
has a minimum $\varepsilon>0$. Now take an iterated barycentric subdivision
$sd^t(K)$ of $K$ with the standard geometric realization
$||sd^t(K)||=||K||$. As the dimension of $sd^t(K)$ is the same as that of $K$,
we can simply consider $sd^t(K)$ with the same map $g:\mathbb
S^h\to||sd^t(K)||$. If 
$t$ is high enough the maximum diameter of the body of a simplex in $sd^t(K)$
is below $\varepsilon/2$ and therefore our assumption on antipodal points is
satisfied.

Let the vertices of $K$ be $[m]=1,\dots,m$. We
define $A_i\subseteq \mathbb S^h$ for all $i$ in $[m]$ by letting ${\mbf x}\in
A_i$ 
if and only if $\alpha_i=\max_j\alpha_j$ in the formal convex combination
$g({\mbf x})=\sum_{j=1}^m\alpha_j||j||$ with $\{j:\alpha_j>0\}\in K$. Clearly,
the closed sets $A_i$ cover $\mathbb S^h$. As ${\mbf x}\in A_i$ implies that
$i$ is a 
vertex of the minimal simplex containing $g({\mbf x})$ the point $\mbf x$ is
contained in at most $l$ of the sets $A_i$, and by our assumption above no set
$A_i$ contains antipodal pairs of points. 
\hfill$\Box$

\medskip

The following corollary is just a restatement of the implication
(iv)$\Rightarrow$(i) of the above lemma for later reference. 

\begin{cor} \label{Qcoind}
For any finite graph $G$ we have $Q(\coind(H(G)))\leq\psi(G)-1$.
\hfill$\Box$
\end{cor}

\medskip
\par
\noindent
{\it Remark 2.}:
Using the fact that any $d$-dimensional simplicial complex
has a  geometric realization in $\mathbb R^{2d+1}$ (cf.\ Theorem~1.6.1 in
\cite{Mat}) and Lemma~\ref{lem:equiB} one can show that
$Q(h)\ge\left\lceil{h\over 2}\right\rceil+1$, i.e., the same lower bound that
Ky Fan's theorem implied in \cite{ST}. Indeed, by Lemma~\ref{lem:equiB}
$Q(h)\leq\left\lfloor{h-1\over 2}\right\rfloor+1$ would imply the existence of
a continuous map $g:\mathbb S^h\to ||K||$ where $K$ is an at most
$\left\lfloor{h-1\over 2}\right\rfloor$-dimensional simplicial complex and
$g(\mbf x)\neq g(-\mbf x)$ for any $\mbf x\in \mathbb S^h$. But $K$ can be
realized in $\mathbb R^h$, so this way we would obtain a continuous map from
$\S^h$ to $\mathbb R^h$ with no coinciding images of antipodal points. This
would contradict the Borsuk-Ulam theorem. 

Using the fact that a $1$-dimensional simplicial complex always have a
universal cover (it is an infinite tree) which can be embedded into $\mathbb
R^2$, the above argument can be extended to prove Shkliarsky's result
\cite{Sklj} stating $Q(2)\ge 3$. This method, however, fails to show $Q(2r)>
r+1$ for $r>1$, which is the most difficult statement in the lower bound part
of Theorem~\ref{qfalse}, cf. \cite{Scep}, \cite{IJ}--\cite{Jaworgen}.   
\hfill$\Diamond$

\medskip

Note that the lower bound $\lceil{h\over 2}\rceil+1\leq Q(h)$
implied by Ky Fan's theorem (cf. \cite{ST}) together with
Corollary~\ref{Qcoind} readily implies a weaker version of
Theorem~\ref{thm:lowb}. Namely, they imply that if $G$ is strongly
topologically $t$-chromatic for some $t\ge 2$, then $\psi(G)\ge\lceil
t/2\rceil+1$.

\medskip
\par
\noindent
{\bf Proof of Theorem~\ref{thm:ps4}.}
By Corollary~\ref{Qcoind} we have $Q(\coind(H(G)))\leq \psi(G)-1$. Using
Theorem~\ref{qfalse} this implies $\psi(G)\ge \left\lfloor t\over
  2\right\rfloor+2$ if $\coind(H(G))=t-2\ge 1$. 
\hfill$\Box$

\medskip
\par
\noindent

Thus any $t$-chromatic Kneser graph, Schrijver graph, generalized Mycielski
graph, or Borsuk graph has local chromatic number at least
$\lfloor t/2\rfloor+2$. For Borsuk graphs it follows immediately from
Lemma~\ref{lem:equiB} that this bound is sharp for $B(t-1,\alpha)$ if
$\alpha<2$ is large enough. 
The results in \cite{ST} imply that it is also 
sharp for many Schrijver graphs and generalized Mycielski graphs. 
This was shown there for odd $t$, while for even $t$ a
gap of $1$ remained in \cite{ST} (compared to the lower bound proven there). 
This gap is closed now. 
Thus we can formulate the
following corollary generalizing Theorems~3 and 5 of
\cite{ST} for the even chromatic case. For the precise meaning of the phrase
``defining parameters'' in the statement below we refer the reader to the
corresponding cited statements of \cite{ST}.

\begin{cor}
Let $t$ be fixed. 
If $G$ is a $t$-chromatic Schrijver graph or a $t$-chromatic generalized
Mycielski graph with large enough defining parameters, then  
$$\psi(G)= \left\lfloor t\over 2\right\rfloor+2.$$
\end{cor}

\proof
The lower bound follows from Theorem~\ref{thm:ps4} and the fact that these
graphs are strongly topologically $t$-chromatic. The matching upper
bound follows from Theorems~3 and 5 in \cite{ST}. \hfill$\Box$

The upper bound is trivial when $t=4$, thus there we have unconditionally,
that any $4$-chromatic Kneser graph, Schrijver graph, generalized
Mycielski graph, or Borsuk graph has local chromatic number $4$.

\medskip
\par
\noindent
{\em Remark 3.}
Let the graph $G$ be a {\em quadrangulation} of a compact two dimensional
surface $R$, i.e., $G$ is drawn in the surface with all the resulting cells
being quadrangles. In this case $H(G)$ is closely related to $R$. In
particular it is easy to show that $\coind(H(G))\ge2$ if $G$ is a
quadrangulation of the projective plane and $G$ is not bipartite. Using
Theorem~\ref{thm:ps4} this implies that the local chromatic number of $G$ is
at least $4$ generalizing the lower bound part of Youngs' result \cite{Y}
which states that such graphs are $4$-chromatic. 
It has been widely studied when quadrangulations of surfaces 
have (ordinary) chromatic number at least $4$, see \cite{AHNNO,MS,Y}. In such
cases four distinct colors can always be found {\em locally}: any proper
coloring has a multicolored quadrangular cell. 
(We call a set of vertices or a subgraph multicolored if every vertex in it
receives a different color.)
Thinking of this four-cycle
as a complete bipartite graph there is a clear connection to what is called
the Zig-zag Theorem in \cite{ST} (cf.\ also Ky Fan's paper \cite{kyfan2}).  
Proving that the local chromatic number is at least $4$ constitutes
finding a different multicolored subgraph: a star with four vertices. This
seems to be harder. The observation that non-bipartite quadrangulations of the
projective plane have local chromatic number at least $4$ generalizes 
to certain quadrangulations of the Klein bottle. Surprisingly, there are
quadrangulations of other surfaces for which a multicolored cell can be found
in every proper coloring but the local chromatic number is only $3$. See the
forthcoming paper \cite{MST} on quadrangulations of surfaces.
\hfill$\Diamond$  

\medskip

In view of the results in \cite{ST} and \cite{ST3} it seems
natural to ask what complete bipartite graphs $K_{k,l}$ must have a
multicolored copy in every proper coloring of any (strongly)
topologically $t$-chromatic graph. To avoid trivialities we always assume
$k,l\ge1$ when speaking about $K_{k,l}$. Using the results in \cite{ST} and
this paper we can give a complete answer for topologically
$t$-chromatic graphs and an almost complete answer for strongly
topologically $t$-chromatic graphs. Note that \cite{ST3} treated the
same problem for proper $t$-colorings of topologically $t$-chromatic
graphs and found a different characterization.

Let us consider topologically $t$-chromatic graphs first. As some of
these graphs are indeed $t$-chromatic we must have $k+l\le t$. By
Corollary~\ref{cor:pelda} (see below in Section~\ref{versus}) the local
chromatic number of some of them is $\lceil t/2\rceil+1$, so we must
also have $k,l\le\lceil t/2\rceil$. For the
remaining graphs $K_{k,l}$ the Zig-zag Theorem of \cite{ST} provides a
positive answer: any proper coloring of a topologically $t$-chromatic
graph contains a multicolored copy of $K_{\lfloor t/2\rfloor,\lceil
t/2\rceil}$ and thus also of its subgraphs. 

A $t$-coloring of a topologically $t$-chromatic graph cannot avoid a
multicolored copy of $K_{k,l}$ for any pair of natural numbers $k,l$ with
$k+l\leq t$ \cite{ST3}, but some $(t+1)$-colorings simultaneously avoid
multicolored copies of all graphs $K_{k,l}\nsubseteq K_{\lfloor t/2\rfloor,
\lceil t/2 \rceil}$. For even $t$ this follows from
Corollary~\ref{cor:pelda}, while for odd $t$ this is stated in \cite{ST}.
(Notice that in this latter case we need to avoid multicolored $K_{\lceil
t/2\rceil,\lceil t/2\rceil}$ subgraphs which does not follow from attaining
local chromatic number $\lceil t/2\rceil+1$.) 

Some strongly topologically $t$-chromatic graphs are also
$t$-chromatic but as we have seen their lowest possible local
chromatic number is $\lfloor t/2\rfloor+2$ (attained by some
Schrijver, Borsuk and generalized Mycielski graphs, see \cite{ST}).
This means that in order to always find a multicolored copy of
$K_{k,l}$ in a proper coloring of a strongly topologically
$t$-chromatic graph we need $k+l\le t$ and $k,l\le\lfloor
t/2\rfloor+1$. Similarly to the previous case, a single $(t+1)$-coloring of a
strongly $t$-chromatic graph can avoid multicolored copies of $K_{k,l}$ for
all $k$, $l$ with $k+l>t$ or $\max(k,l)>\lfloor t/2\rfloor+1$. For $t$ odd this
is proven in \cite{ST}. It easily extends
also to the even $t$ case by taking the Mycielskian $M(G)$ of 
a strongly topologically $(t-1)$-chromatic graph $G$ with a coloring of the
above type and extend this to a proper coloring of $M(G)$ (which is a strongly
topologically $t$-chromatic graph) in the following way. We keep the original
coloring in the first layer, introduce a new color for all the second layer
and use one of the old colors for the top vertex. (For the definition of
Mycielskians we refer again to \cite{SchU} or \cite{ST}.)
 
For most of the remaining complete bipartite graphs the Zig-zag theorem implies
the existence of a multicolored version, the only case not covered is
that of $K_{t/2+1,l}$ for even $t$ and $1\le l\le t/2-1$.
Theorem~\ref{thm:ps4} is equivalent to an affirmative answer in the
$l=1$ case. For $l>1$ we do not know the answer. Here we ask the
problem in the strongest possible form corresponding to $l=t/2-1$
\medskip
\par
\noindent
{\em Question.} Let $t\ge 6$ be an even integer. Is it true
that if a strongly topologically $t$-chromatic graph is properly colored (with
any number of colors) then it always contains a multicolored
$K_{t/2-1,t/2+1}$ subgraph? 
\medskip

Using similar techniques to those used in \cite{ST} and \cite{ST3} and also
in this paper, an
affirmative answer would immediately follow from an affirmative answer to the
following topological analog of the above question. 

\medskip
\par
\noindent
{\em Topological question.} Let $h\ge 4$ be even and let the sphere $\mathbb
S^h$ be covered by open sets $A_1,\dots,A_m$ that satisfy $A_i\cap
(-A_i)=\emptyset$ for all $i$. Is it true that there always exists an $x\in
\mathbb S^h$ such that $x$ is covered by at least $h/2+2$ and $-x$ is covered
by at least $h/2$ different $A_i$'s?

\medskip
Note that an affirmative answer would give a strengthening of the lower bound
part of Theorem~\ref{qfalse}, while the two are equivalent if we set $h=2$.

\section{Topological $t$-chromaticity versus strong topological
$t$-chromaticity} \label{versus} 

In this section we compare topological $t$-chromaticity and strong topological
$t$-chromaticity, especially in their implications to the
local chromatic number. 

As stated in (\ref{eq:chib1}) strong topological $t$-chromaticity
implies topological $t$-chromaticity, which, in turn, implies that the
graph is indeed at least $t$-chromatic. It is easy to see that for $t=2$ or
$3$ both topological conditions are equivalent with the graph
having chromatic number at least $t$. 
This is not the case for $t\ge 4$ as it follows from an observation by
Walker \cite{Walk} made also by Matou\v{s}ek and Ziegler \cite{MZ}. This
observation (in terms of \cite{MZ}) is that any graph $G$
without a 
$4$-cycle satisfies $\ind(||B(G)||)\leq 1$. Using the already mentioned
$\2$-homotopy equivalence of $H(G)$ and $||B(G)||$ and the result of Erd\H os
\cite{erdos} that there exist graphs with arbitrarily high chromatic number and
girth this shows that the two sides of the first inequality in
(\ref{eq:chib1}) can be arbitrarily far apart. 

If one of the other three inequalities in (\ref{eq:chib1}) is
strict then we have $\ind(H(G))>\coind(H(G))$. $\2$-spaces having
different index and coindex are called {\em nontidy} by Matou\v sek
\cite{Mat}. Constructing such spaces
do not seem obvious but such constructions are known, see, e.g., a list in
\cite{Mat}, page 100. Csorba \cite{Cs} and \v{Z}ivaljevi\'c \cite{Ziv} proved
that for any finite free $\2$-complex $K$ there exists a finite graph $G$ such
that $||B(G)||$ (and thus also $H(G)$) is $\2$-homotopy equivalent to
$||K||$. Some of the nontidy spaces, e.g., the projective space
${\mathbb R}P^{2i-1}$ with a suitable involution, have a triangulation
(i.e., it is $\2$-homeomorphic to the body of a finite $\2$-complex). So we
have 
examples of graphs $G$ with $\ind(H(G))>\coind(H(G))$, and (from the
properties of ${\mathbb R}P^{2i-1}$) even $\coind(H(G))=1$ with $\ind(H(G))$
arbitrarily high. This shows that the difference between the two sides of at
least one of the
second, third, or last inequality of (\ref{eq:chib1}) is unbounded (but, as
mentioned below, it certainly cannot be the second). Further
study of the space $\susp({\mathbb R}P^{2i-1})$ shows that its coindex is $2$
\cite{Zivketlap,Feherszo}, showing that for the above graphs there is an
unbounded difference between the two sides of the third inequality in
(\ref{eq:chib1}), while the last inequality holds with equality.

Based on another example appearing in \cite{Mat}, page 100, constructed by
Csorba, Matou\v{s}ek, and \v{Z}ivaljevi\'c, an example of a $\2$-space $X$ is
demonstrated by Csorba \cite{Cs} which satisfies
$\ind(X)=\ind(\susp(X))$. Since this space can also be triangulated,
it shows the existence of graphs for which the second inequality is strict in
(\ref{eq:chib1}) (using again the above mentioned result of
Csorba \cite{Cs} and \v{Z}ivaljevi\'c \cite{Ziv}). Nevertheless, as $H(G)$ is
contained in $B_0(G)$ the sides of the second inequality can differ by at
most $1$.

Our main concern is the last inequality of (\ref{eq:chib1})
which is 
between the defining quantities of topological and strongly topological
$t$-chromaticity. Here we show not only the possibility of strict inequality,
but also the existence of a topologically $t$-chromatic but not strongly
topologically $t$-chromatic graph for which
Theorem~\ref{thm:lowb} is tight while $t$ is even. (For odd $t$ several 
examples are shown in \cite{ST} for the tightness of the lower bound in
Theorem~\ref{thm:lowb}, however, those examples are also strongly
topologically $t$-chromatic.) In case of $t=2r$ this means that our graph 
has local chromatic number $r+1$ in contrast to strongly topologically
$2r$-chromatic graphs for which the local chromatic number must be at least
$r+2$ 
according to Theorem~\ref{thm:ps4}. Thus our examples will not only separate
topological $2r$-chromaticity from strong topological $2r$-chromaticity
but show 
that the difference is in fact relevant also in terms of its consequences for
the local chromatic number. 
We do not have examples where the sides of the last
inequality of (\ref{eq:chib1}) differ by more than $1$.

\smallskip

Our examples of topologically $t$-chromatic graphs with local chromatic number
equal to $\left\lceil{t\over 2}\right\rceil+1$, the lower bound in
Theorem~\ref{thm:lowb}, are the universal graphs $U(2r-1,r)$ defined below in
the more general setting as they appear in \cite{EFHKRS}. From now on we keep 
using the notation $[m]=\{1,\dots,m\}$ already introduced in the proof of
Lemma~\ref{lem:equiB}. 

\begin{defi}\label{defi:Umr} {\rm (\cite{EFHKRS})}
For positive integers $r\leq m$ we define the graph $U(m,r)$ as follows.
\begin{eqnarray*}
V(U(m,r))&=&\{(i,A): i\in [m], A\subseteq [m], |A|=r-1, i\notin A\}\\
E(U(m,r))&=&\{\{(i,A),(j,B)\}: i\in B, j\in A\}
\end{eqnarray*}
\end{defi}

The graphs $U(m,r)$ characterize local chromaticity in the sense that a graph
$G$ satisfies $\psi(G)\leq r$, and this value can be attained by a coloring
with at most $m$ colors, if and only if there is a homomorphism from $G$ to
$U(m,r)$ (see Lemma 1.1 in \cite{EFHKRS}). In
particular, it is easy to find the coloring showing $\psi(U(m,r))\leq r$: for
each vertex $(i,A)$ use $i$ as its color. We refer to this coloring as the
natural coloring of $U(m,r)$.
(Note that $\chi(U(m,r))<m$ whenever $m>r$, cf. \cite{EFHKRS}, thus this is
not an optimal coloring concerning the number of colors used. In fact, it is
easy to see that if $\psi(G)<\chi(G)$ then any coloring of $G$ attaining
$\psi(G)$ must use more than $\chi(G)$ colors. The reason is that in a proper
coloring of $G$ with $\chi(G)$ colors each color class must contain a vertex
which has a neighbor in all other color classes. Otherwise the color class
with no such vertex could be eliminated resulting in a proper coloring with
less than $\chi(G)$ colors.) 

\medskip
\par
\noindent
{\em Remark 4.}
The above discussion shows that the local chromatic number fits into the
framework described in Chapter 1 of Kozlov's survey \cite{Kozl}. Namely, that 
$\psi(G)$ could also be defined as the minimum $r$ for which $G$ admits a
homomorphism into one of the graphs $U(m,r)$. In the language of \cite{Kozl}
this defines $\psi(G)$ via the state graphs $U(m,r)$ and valuation
$U(m,r)\mapsto r$. 
\hfill$\Diamond$

\medskip

To be able to speak about topological $t$-chromaticity with respect to the
graph $U(m,r)$ we need to consider $B_0(U(m,r))$. It is going to be useful to
introduce an exponentially smaller $\2$-equivalent complex.

\begin{defi} \label{defi:Lmr}
Let $K_m$ denote the complete graph on the vertex set $V(K_m)=[m]$ and let
$L_m=B_0(K_m)$. For a positive integer $r\le m$ let
$L_{m,r}$ denote the subcomplex of $L_m$ that consists of those
simplices $S\uplus T$, for which $|S|<r$ and $|T|<r$. Let
$L_{m,r}'=L_{m,r}\cup\{S\uplus\emptyset:S\subseteq[m]\}\cup\{\emptyset\uplus
T:T\subseteq[m]\}$. The bodies $||L_{m,r}||$ and $||L_{m,r}'||$ are
$\2$-spaces with the involution inherited from $||L_m||$.
\end{defi}

\begin{lem} \label{lem:uaz}
For every $m$ and $r$ we have $||L_{m,r}'||\leftrightarrow||B_0(U(m,r))||$.
\end{lem}

\proof
We have a simplicial $\2$-map $B_0(U(m,r))\to L_{m,r}'$ given by
$+(i,A)\mapsto+i$ and $-(i,A)\mapsto-i$. This shows
$||B_0(U(m,r))||\to||L_{m,r}'||$.

We give a monotonously decreasing map $g$ from the
simplices in the barycentric subdivision $sd(L_{m,r}')$ to the
simplices of $B_0(U(m,r))$. This map can be considered as a
simplicial map from the second subdivision $sd(sd(L_{m,r}'))$ to the
subdivision $sd(B_0(U(m,r)))$ and thus $||g||$ 
(the piecewise linear extension of $g$) maps
$||sd(sd(L_{m,r}'))||=||L_{m,r}'||$ to
$||sd(B_0(U(m,r)))||=||B_0(U(m,r))||$. This is clearly a $\2$-map showing
$||L_{m,r}'||\to||B_0(U(m,r))||$ as stated.

Recall that the vertices
of $sd(L_{m,r}')$ are the simplices of $L_{m,r}'$ and a non-empty set of
vertices forms a simplex in $sd(L_{m,r}')$ if it is linearly ordered by
inclusion. Let therefore $C$ be a simplex of $sd(L_{m,r}')$ and let
$S\uplus T$  be its smallest vertex and $S'\uplus T'$ be its largest
vertex. We set $g(C)=W\uplus Z$ with $W=\{(i,H)\in
V(U(m,r)):i\in S, T'\subseteq H\}$ and $Z=\{(i,H)\in V(U(m,r)):i\in T,
S'\subseteq H\}$. Any
pair of vertices $w\in W$ and $z\in Z$ are connected in $U(m,r)$, so
$g(C)\in B_0(U(m,r))$. The map $g$ is clearly monotonously
decreasing. Thus simplices of $sd(sd(L_{m,r}'))$ are mapped into simplices of
$sd(B_0(U(m,r)))$ provided $g(C)$ is not empty. Assume first
that $S\ne\emptyset$. We have $S\subseteq S'\ne\emptyset$, so by the
definition of $L_{m,r}$ we have $|T'|\le r-1$. We choose $i\in S$ and a set
$H\supseteq T'$ with $|H|=r-1$ and $i\notin H$. We have $(i,H)\in
W$, so $W\ne\emptyset$. The same argument shows that if $T\ne\emptyset$ then
$Z\ne\emptyset$. As $S\uplus T$ is a simplex either $S\ne\emptyset$ or
$T\ne\emptyset$, and we have $g(C)\ne\emptyset$ in either case.
\hfill$\Box$
\medskip

\medskip
\par
\noindent
{\it Remark 5.} Though we need only the above proven $\2$-equivalence of
$||L_{m,r}'||$ and $||B_0(U(m,r))||$, we mention that they are actually
$\2$-homotopy equivalent. 

To prove this we show that the $\2$-mappings
$||f||:||B_0(U(m,r))||\to||L'_{m,r}||$ and
$||g||:||L'_{m,r}||\to||B_0(U(m,r))||$ satisfy that both
$||f||\circ ||g||$ and $||g||\circ ||f||$ are $\2$-homotopic to the
identity of the respective spaces. Here $\2$-homotopic means that
they are homotopic with every layer of the homotopy being a $\2$-map and 
$f$ is the simplicial map corresponding to the natural coloring of $U(m,r)$,
while $g$ denotes the same map as in the proof above. 

Let $a$ and $b$ be $\2$-maps from $X$ to $||C||$ where $C$ is a simplicial
$\2$-complex. If for all $x\in X$ the points
$a(x)$ and $b(x)$ are contained in the body of a common simplex, then linear
interpolation between $a$ and $b$ proves that
they are $\2$-homotopic.

This simple observation can be directly used to show that
$||f||\circ||g||$ is $\2$-homotopic with the identity on
$||L'_{m,r}||$. Unfortunately, the same argument cannot be used
directly to show that $||g||\circ||f||$ is homotopic to the identity
on $||B_0(U(m,r))||$. We introduce the simplicial map
$h:sd(sd(B_0(U(m,r))))\to sd(B_0(U(m,r)))$ mapping a chain of
simplices from $B_0(U(m,r))$ to its smallest element. Clearly,
$||h||:||B_0(U(m,r))||\to||B_0(U(m,r))||$ is a $\2$-map. Now the elementary
argument can be used to show that both the identity and
$||g||\circ||f||$ are $\2$-homotopic to $||h||$. This shows that they
are $\2$-homotopic to each other, too.
\hfill$\Diamond$

\medskip
In the following lemma we use the notion of Bier spheres. For a complex $K$
with $V(K)\subseteq[m]$, $[m]\notin K$ its Bier sphere is defined as
$$\hbox{Bier}_m(K)=\{S\uplus T\in L_m:S\in K,\overline T\notin K\},$$ 
where $\overline T=[m]\setminus T$ is the complement of $T$. The basic result
on Bier spheres is that they are always triangulations of a sphere:
$||\hbox{Bier}_m(K)||\cong \mathbb S^{m-2}$. For a proof of this result see,
e.g., 
Theorem~5.6.2 in \cite{Mat}, or \cite{deLong}. 

\begin{lem} \label{lem:Bier}
For $r\ge1$ we have $||L_{2r-1,r}||\cong \mathbb S^{2r-3}$.
\end{lem}

\proof
Observe that $L_{2r-1,r}$ is just the Bier sphere $\hbox{Bier}_{2r-1}(K)$ of
the simplicial complex $K={[2r-1]\choose {\leq r-1}}$ consisting of the at 
most $(r-1)$-element subsets of $[2r-1]$.
\hfill$\Box$

\begin{cor} \label{cor:coindB}
The graph $U(2r-1,r)$ is topologically $(2r-2)$-chromatic.   
In particular we have $$\coind(||B_0(U(2r-1,r))||)=2r-3.$$
\end{cor}

\proof
By Lemma~\ref{lem:uaz} we have
$\coind(||B_0(U(2r-1,r))||)=\coind(||L'_{2r-1,r}||)$. By containment we have
$\coind(||L'_{2r-1,r}||)\ge\coind(||L_{2r-1,r}||)$. By Lemma~\ref{lem:Bier}
(and the remark in the introductory part of Subsection~\ref{topre} about
homotopy spheres) we have $\coind(||L_{2r-1,r}||)=2r-3$. 

The reverse inequality follows from applying the inequality 
$\chi(G)\ge \coind(||B_0(G)||)+1$ to $G=U(2r-1,r)$ and using the inequality
$\chi(U(2r-1,r))\leq 2r-2$. The latter is a special case of the fact mentioned
above, that $\chi(U(m,r))<m$ if $r<m$. 
\hfill$\Box$

\medskip
\par
\noindent
{\em Remark 6.}
The fact that $\chi(U(2r-1,r))\ge 2r-2$ is a special case of Theorem~2.6 in
\cite{EFHKRS}. 
This remark parallels Remark~3 in \cite{ST} which explains how the upper bound
results of \cite{ST} imply $\chi(U(2r,r+1))\ge 2r-1$, another special case of
Theorem~2.6 in \cite{EFHKRS}. In \cite{ST} this follows from the proof of local
$(r+1)$-chromaticity of some strongly topologically $(2r-1)$-chromatic graphs
that can be attained by using $2r$ colors. This implies the existence of
homomorphisms from some strongly topologically $(2r-1)$-chromatic graphs to
$U(2r,r+1)$. Besides implying $\chi(U(2r,r+1))\ge 2r-1$ this also shows that
the graphs $U(2r,r+1)$ are strongly topologically $(2r-1)$-chromatic. (Their
chromatic number is $2r-1$, indeed, by the same argument we have in the second
part of the proof of Corollary~\ref{cor:coindB}.) 
The above is in contrast to the case of $U(2r-1,r)$, since these graphs, as
stated below in Corollary~\ref{cor:pelda}, are only topologically
$(2r-2)$-chromatic but not strongly topologically $(2r-2)$-chromatic. 
\hfill$\Diamond$

\begin{cor} \label{cor:pelda}
For any $r\ge1$ there exists a topologically $2r$-chromatic graph with
local chromatic number $\psi(G)=r+1$ which can be attained by a
$(2r+1)$-coloring. In particular, topological $2r$-chromaticity implies
neither strong topological $2r$-chromaticity nor 
that the local chromatic number is at least $r+2$.
\end{cor}

\proof The example claimed is $U(2r+1,r+1)$. The local chromatic number
is attained by its natural coloring. Topological $2r$-chromaticity is given
by Corollary~\ref{cor:coindB}. 
Theorem~\ref{thm:ps4} shows that $U(2r+1,r+1)$ is not strongly
topologically $2r$-chromatic, since then its local chromatic number should
be larger.  
\hfill$\Box$

\section{Direct separation arguments}
\label{remark14}

Inequality~(\ref{eq:chib1}) and our statement that $U(2r+1,r+1)$
satisfies topological $2r$-chromaticity, but not strong topological
$2r$-chromaticity shows that $H(U(2r+1,r+1))$ has different index and
coindex. While, as we already mentioned, the existence of such
spaces has been known (even with arbitrarily high difference between the index
and 
the coindex, see page 100 of \cite{Mat} and the references therein),
$H(U(5,3))$ 
yields a particularly simple and elementary example. See the argument below on
compact orientable $2$-manifolds. 

First we claim a variant of Lemma~\ref{lem:uaz} for the hom space. Let
$H_{m,r}=||L_{m,r}||\cap H(K_m)$. We claim that
$H(U(m,r))\leftrightarrow H_{m,r}$. The proof is
almost identical to that of Lemma~\ref{lem:uaz}.

Notice that $H_{2r+1,r+1}$ is a topological $(2r-2)$-manifold. To see this
consider $H_{m,r}$ as the body of the cell complex $\hat
H_{m,r}={\rm Hom}(K_2,K_m)\cap L_{m,r}$. It is enough to verify that it is
connected 
and the {\em link} of any vertex is a triangulation of the same sphere. Here
the link of a vertex $V$ in the complex $K$
consists of the sets $W\setminus V$ for cells $W$ in $K$ containing $V$. Note
that the link of a vertex in the cell complex $\hat H_{m,r}$ is a
simplicial complex. By the symmetry of $\hat H_{m,r}$ the link of each vertex 
is isomorphic. The link of the vertex $\{m-1\}\uplus\{m\}$ is
$L_{m-2,r-1}$. In case of $\hat H_{2r+1,r+1}$ this link is $L_{2r-1,r}$ and
thus it is a triangulated $\mathbb S^{2r-3}$ by Lemma~\ref{lem:Bier} as needed.

As $H_{2r+1,r+1}$ is a $(2r-2)$-manifold embedded in
$||L_{2r+1,r+1}||\cong {\mathbb S}^{2r-1}$ it is orientable. One can 
easily compute the Euler-characteristic of $H_{5,3}$ directly: its defining
cell complex has $20$ vertices, $30$ cells of dimension $2$, and $60$ edges,
thus the Euler-characteristic is $-10$.  
This shows that $H_{5,3}$ is the orientable compact $2$-manifold of
genus $6$. Consider this manifold
as a sphere $\mathbb S^2$ with six ``handles'' arranged in a centrally
symmetric 
manner. The central reflection gives the involution of the space.

The following argument is a direct and simple proof that any compact
orientable $2$-manifold $T$ of even and positive genus and with
the involution as above satisfies that its index is $2$, while its coindex is
$1$. We will show that the index is at least $2$ by showing that the coindex
of its suspension is at least $3$. 
Thus by the result of Csorba \cite{Cs} and \v{Z}ivaljevi\'c \cite{Ziv} and
since these spaces admit triangulations, each of these examples yield
different topologically $4$-chromatic graphs that are not strongly
topologically $4$-chromatic.

Note that $\ind(T)\le2$ and $\coind(T)\ge1$ are trivial.

We prove $\coind(\susp(T))\ge3$ by giving the
explicit mapping. Then $\ind(\susp(T))\ge3$ by the Borsuk-Ulam Theorem and
$\ind(T)\ge2$ follows from Lemma~\ref{suspension}. We choose $T$ as a
subspace of $\mathbb S^3$ closed for the
involution. This can be done in a smooth way such that the points ${\mbf x}\in
\mathbb S^3$ within some distance $\varepsilon>0$ to $T$ has a unique closest
point 
$\hat{\mbf x}\in T$. We denote by $T^+$ and $T^-$ the two components of
$\mathbb S^3\setminus T$. If we identify the points in $T^+$ far away from $T$
and also identify the points in $T^-$ far away from $T$, then the resulting
factor space is naturally homeomorphic to $\susp(T)$. The resulting map
$f:\mathbb S^3\to\susp(T)$ can be given as follows.
$$f({\mbf x})=\left\{\begin{array}{lll}
(*,1)&&\hbox{if }\hbox{dist}({\mbf x},T)\ge\varepsilon,{\mbf x}\in T^+\\
(\hat{\mbf x},\hbox{dist}({\mbf x},T)/\varepsilon)&&
\hbox{if }\hbox{dist}({\mbf x},T)<\varepsilon,{\mbf x}\in T^+\\
({\mbf x},0)&&\hbox{if }{\mbf x}\in T\\
(\hat{\mbf x},-\hbox{dist}({\mbf x},T)/\varepsilon)&&
\hbox{if }\hbox{dist}({\mbf x},T)<\varepsilon,{\mbf x}\in T^-\\
(*,-1)&&\hbox{if }\hbox{dist}({\mbf x},T)\ge\varepsilon,{\mbf x}\in T^-
\end{array}\right.$$

We prove $\coind(T)<2$ by a similar argument as the one hinted in Remark 2. 
A continuous map
$f:\mathbb S^2\to T$ lifts to the universal covering space ${\mathbb R}^2$ of
$T$, but 
by the Borsuk-Ulam Theorem $\hat f:\mathbb S^2\to{\mathbb R}^2$ identifies two
antipodal points of the sphere, so $f$ also identifies two antipodal points,
and therefore $f$ is not a $\2$-map.

It is worth noting where the argument showing $\coind(\susp(T))\ge 3$ fails
for compact orientable $2$-manifolds $T$ of odd genus (with the involution
given by the reflection in their standard self-dual embeddings in $\mathbb
S^3$). The 
function $f$ defined as above is
{\em not} a $\2$-map because the involution on $\mathbb S^3$ does not switch
the 
components of $\mathbb S^3\setminus T$ in this case. Indeed, we have
$\ind(T)=1$ for 
such surfaces $T$.

The argument above shows that for the manifold $T$ considered above one has
$\coind(\susp(T))-\coind(T)=2$. It would be interesting to find spaces $T$
with $\coind(\susp(T))-\coind(T)$ arbitrarily large.

\medskip

The space $T=H_{2r+1,r+1}$ provides an example of a $\2$-space $T$ with
$\coind(T)\le 2r-3$ and $\coind(\susp(T))\ge 2r-1$. The following lemma 
can be used to find examples for spaces $T$ with $\coind(T)\le d-2$ and
$\coind(\susp(T))\ge d$ also for even values of $d$. 

Below we use cohomologies over $\2$. 

\begin{lem}\label{cohom}
If a compact $d$-manifold $T$ has a non-trivial cohomology $l$ in some
dimension $1\le i\le d-1$ and it is the body of a free simplicial $\2$-complex,
then $\coind(T)<d$.
\end{lem}

\proof
We need to show that no $\2$-map $f:\mathbb S^d\to T$ exists.

Assume for a contradiction that such a map $f$ exists. It induces a reverse map
$f^*$ on the cohomologies, in particular it maps $l$ to an
$i$-cohomology of $\mathbb S^d$. As no such non-trivial cohomology exists, we
have 
$f^*(l)=0$. By Poincar\'e duality there exists a $(d-i)$-cohomology $l'$
with the cup product $l\!\smile\!l'=z$ being the only non-trivial
$d$-cohomology in $T$. As $f^*$ preserves the cup product we have $f^*(z)=0$.

As $T$ is the body of a $d$-dimensional free simplicial $\2$-complex there
is a $\2$-map $g:T\to \mathbb S^d$. Let $w$ be the only non-trivial
$d$-cohomology of 
$\mathbb S^d$. The homomorphism $g^*$ induced by $g$ maps $w$ either to $z$ or
to $0$, 
in either case $(g\circ f)^*(w)=f^*(g^*(w))=0$. This shows
that $g\circ f:\mathbb S^d\to \mathbb S^d$ has even {\em degree} which
contradicts the fact 
that it is a $\2$-map (cf. \cite{Bred}, Theorem 20.6 on page 244). 
This contradiction proves $\coind(T)<d$. 
\hfill$\Box$

\medskip

Let $T$ be a $(d-1)$-manifold obtained by attaching two homeomorphic
``handles'' to the sphere $\mathbb S^{d-1}$.
We start with $\mathbb S^{d-1}\subset \mathbb S^d$ and
attach the handles inside $\mathbb S^d$ in a 
centrally symmetric way and smoothly, just as in the $d=3$ case earlier. The
central reflection of $\mathbb S^d$ gives the involution in $T$. We can prove
$\coind(\susp(T))\ge d$ via the same explicit $\2$-map $f:\mathbb
S^d\to\susp(T)$ as in the $d=3$ case.
Note that the space $T$ constructed admits a triangulation where the
involution is simplicial. 
As $T$ is a $(d-1)$-manifold with the above property, by
Lemma~\ref{cohom} it is enough to find a nontrivial cohomology (over $\2$) of
$T$ in some dimension 
$1\le i\le d-2$ and this implies $\coind(T)\le d-2$. By choosing for example
handles homeomorphic to a punctured ${\mathbb S}^k\times {\mathbb S}^l$ with
$k,l\ge1$, $k+l=d-1$ one makes sure that the $k$th and $l$th cohomology groups
are non-trivial. 
More formally, we could say that we construct the manifold $T$
from the sphere $\mathbb S^{d-1}$ by applying two $(k-1)$-surgeries
($1\leq k \leq d-2$), in a centrally symmetric way.

\medskip

For every even $t>2$ we have found a graph that is topologically $t$-chromatic
but not strongly topologically $t$-chromatic: it is $U(t+1,t/2+1)$. As
we pointed out in Remark~6, the analogous graphs $U(t+1,(t+1)/2+1)$ with $t$
odd are strongly topologically $t$-chromatic and their chromatic number is
also $t$. 
Still, using the $\2$-space $T$ constructed in the last paragraph with
$\coind(T)\le d-2$ and $\coind(\susp(T))\ge d$ one can separate
topological $t$-chromaticity from strong topological $t$-chromaticity in a
similar way also for odd $t$. Choose $d=t-1$. Using again that the space $T$
constructed admits a triangulation where the involution is simplicial we get
by the 
results of Csorba \cite{Cs} and \v{Z}ivaljevi\'c \cite{Ziv} that there exists a
finite graph $G$ with $H(G)$ $\2$-(homotopy) equivalent to $T$ and therefore
$||B_0(G)||$ $\2$-(homotopy) equivalent to $\susp(T)$. This graph $G$ is
topologically $t$-chromatic, but not strongly topologically $t$-chromatic. 
Thus we proved
\begin{cor}
For every integer $t>2$ there exists a graph that is
topologically $t$-chromatic, but not strongly topologically $t$-chromatic.  
\hfill$\Box$
\end{cor}

\bigskip
\par
\noindent
{\bf Acknowledgments:}
We are grateful to Imre B\'ar\'any, P\'eter Csorba, G\'abor Elek, L\'aszl\'o
Feh\'er, L\'aszl\'o Lov\'asz, Ji\v{r}\'{\i} Matou\v{s}ek, G\'abor
Moussong, Andr\'as Sz\H ucs, and Rade \v{Z}ivaljevi\'c 
for many clarifying conversations and e-mail messages that improved our
understanding of the topological concepts used in this paper.


\begin{thebibliography}{99}

\bibitem{AF} J. M. Aarts, R. J. Fokkink, Coincidence and the colouring of
  maps, {\em Bull. London Math. Soc.} {\bf 30} (1998),  no. 1, 73--79.  

\bibitem{AFL} N. Alon, P. Frankl, L. Lov\'asz, 
The chromatic number of Kneser hypergraphs, {\em Trans.\ Amer.\
Math.\ Soc.}, {\bf 298} (1986), 359--370.

\bibitem{AHNNO} D. Archdeacon, J. Hutchinson, A. Nakamoto, S. Negami, K. Ota,  
Chromatic numbers of quadrangulations on closed surfaces, 
{\em J. Graph Theory}, {\bf 37} (2001), no.\ 2, 100--114.

\bibitem{BK} E. Babson, D.\,N. Kozlov, Complexes of graph homomorphisms,   
{\em Israel J. Math.}, {\bf 152} (2006), 285--312, arXiv:math.CO/0310056.  

\bibitem{Bar} I. B\'ar\'any, A short proof of Kneser's
conjecture {\em J. Combin.\ Theory Ser.\ A}, {\bf 25} (1978), no. 3,
325--326.

\bibitem{b} I. B\'ar\'any, personal communication.

\bibitem{Bjhand} A. Bj\"orner, Topological methods, in: {\em Handbook of
Combinatorics} (Graham, Gr\"otschel, Lov\'asz eds.), 1819--1872, Elsevier,
Amsterdam, 1995. 

\bibitem{BjLo} A. Bj\"orner, M. de Longueville, Neighborhood complexes of
stable Kneser graphs, {\em Combinatorica}, {\bf 23} (2003), no.\ 1, 23--34.  

\bibitem{Bred} G. Bredon, {\em Topology and Geometry}, Graduate Texts in
Mathematics 139, Springer-Verlag, New York, 1993.

\bibitem{Cs} P. Csorba, Homotopy types of box complexes, to appear in {\em
    Combinatorica}, arXiv:math.CO/0406118, June 2004.

\bibitem{deLong} M. de Longueville, Bier spheres and barycentric subdivision,
{\em J. Combin.\ Theory Ser.\ A}, {\bf 105} (2004), 355-357. 

\bibitem{erdos} P. Erd\H{o}s, Graph theory and probability, 
{\em Canad.\ J. Math.}, {\bf11} (1959), 34--38.

\bibitem{EFHKRS} P. Erd\H{o}s, Z. F\"uredi, A. Hajnal, P. Komj\'ath,
V. R\"odl, \'A. Seress, Coloring graphs with locally few colors, {\em Discrete
Math.}, {\bf 59} (1986), 21--34.

\bibitem{kyfan} K. Fan, A generalization of Tucker's combinatorial lemma with
topological applications, {\em Annals of Mathematics}, {\bf56} (1952), no.\ 2,
431--437.

\bibitem{kyfan2} K. Fan, Evenly distributed subsets of $\mathbb S^n$ and a
combinatorial application, {\em Pacific J. Math.}, {\bf98} (1982), no.\ 2,
323--325.

\bibitem{Feherszo} L. Feh\'er, personal communication.

\bibitem{Fur} Z. F\"uredi, Local colorings of graphs (Gr\'afok lok\'alis
sz\'{\i}nez\'esei), manuscript in Hungarian, September 2002.

\bibitem{GR} C. Godsil, G. Royle, {\em Algebraic Graph Theory}, Graduate Texts
in Mathematics 207, Springer-Verlag, New York, 2001.

\bibitem{GyJS} A. Gy\'arf\'as, T. Jensen, M. Stiebitz, On graphs with 
strongly independent colour-classes, {\em J. Graph Theory}, {\bf 46} (2004),
1--14.  

\bibitem{Hat} A. Hatcher, {\em Algebraic Topology}, Cambridge University
Press, 2002. Electronic version available at {\tt
http://www.math.cornell.edu/\~{}hatcher/AT/ATpage.html}. 

\bibitem{IJ} M. Izydorek, J. Jaworowski, Antipodal coincidence for maps of
  spheres into complexes, {\em Proc. Amer. Math. Soc.}, {\bf 123} (1995),
  1947--1950. 

\bibitem{Jaworpre}
J. Jaworowski, Existence of antipodal coincidence for maps of spheres,
preprint. 

\bibitem{Jaworgen}
J. Jaworowski, Periodic coincidence for maps of spheres,
{\em Kobe J. Math.}, {\bf 17} (2000), no. 1, 21--26. 

\bibitem{GK} G. Kalai, personal communication.

\bibitem{Kozl} D.\ N.\ Kozlov, Chromatic numbers, morphism complexes, and
  Stiefel-Whitney characteristic classes. {\em Geometric Combinatorics},
  IAS/Park City Mathematics Series 14, American Mathematical Society,
  Providence, RI; Institute for Advanced Study (IAS), Princeton, NJ, to
  appear, arXiv:math.AT/0505563.

\bibitem{KPS} J. K\"orner, C. Pilotto, G. Simonyi, Local chromatic number and
Sperner capacity, {\em J. Combin. Theory, Ser B.}, {\bf 95} (2005), 101--117.  

\bibitem{Kriz}I. K\v r\'{\i}\v z, Equivariant cohomology and lower bounds for
chromatic numbers, {\em Trans.\ Amer.\ Math.\ Soc.}, {\bf 333} (1992), no.\ 2,
567--577; I. K\v{r}\'{\i}\v{z}, A correction to: "Equivariant
cohomology and lower bounds for chromatic numbers", {\em
Trans.\ Amer.\ Math.\ Soc.}, {\bf 352} (2000), no.\ 4, 1951--1952.  

\bibitem{LLKn} L. Lov\'asz, Kneser's conjecture, chromatic number, and
homotopy, {\em J. Combin.\ Theory Ser.\ A}, {\bf 25} (1978), no.\ 3, 319--324. 

\bibitem{LLgomb} L. Lov\'asz, Self-dual polytopes and the chromatic number of
distance graphs on the sphere, {\em Acta Sci.\ Math.\ (Szeged)}, {\bf 45}
(1983), 317--323. 

\bibitem{Mat} J. Matou\v{s}ek, {\em Using the Borsuk-Ulam Theorem. Lectures on
Topological Methods in Combinatorics and Geometry}, Springer-Verlag, Berlin
etc., 2003.  

\bibitem{MZ} J. Matou\v{s}ek, G.\,M. Ziegler, Topological lower bounds for the
chromatic number: A hierarchy, {\em Jahresber. Deutsch. Math.-Verein.}, {\bf
106} (2004), no. 2, 71--90, arXiv:math.CO/0208072.

\bibitem{MS} B. Mohar, P.\,D. Seymour,
Coloring locally bipartite graphs on surfaces, {\em J. Combin.\ Theory Ser.\
B},  {\bf 84} (2002), no.\ 2, 301--310.

\bibitem{MST} B. Mohar, G. Simonyi, G. Tardos, On the local chromatic number
  of quadrangulations of surfaces, manuscript in preparation.

\bibitem{Scep} E. V. \v S\v cepin,
A certain problem of L. A. Tumarkin. (Russian)
{\em Dokl. Akad. Nauk SSSR}, {\bf 217} (1974), 42--43; English translation:
{\em Soviet Math. Dokl.}, {\bf 15} (1974), no. 4, 1024--1026 (1975). 

\bibitem{SchU} E.\,R. Scheinerman, D.\,H. Ullman, {\em Fractional Graph
Theory}, Wiley-Interscience Series in Discrete Mathematics and
Optimization, John Wiley and Sons, Chichester, 1997.

\bibitem{Schr} A. Schrijver, Vertex-critical subgraphs of Kneser graphs, 
{\em Nieuw Arch.\ Wisk.\ (3)}, {\bf 26} (1978), no.\ 3, 454--461.

\bibitem{Sklj} D. Shkliarsky, 
On subdivisions of the two-dimensional sphere. (Russian. English summary)
{\em Rec. Math. [Mat. Sbornik] N. S.} {\bf 16(58)}, (1945), 125--128. 

\bibitem{ST} G. Simonyi, G. Tardos, Local chromatic number, Ky Fan's theorem,
and circular colorings, {\em Combinatorica}, {\bf 26} (2006), 587--626,
arXiv:math.CO/0407075.

\bibitem{ST3} G. Simonyi, G. Tardos, Colorful subgraphs of Kneser-like graphs,
  to appear in {\em European J. Combin.}, arxiv:math.CO/0512019.

\bibitem{STB} G. Simonyi, G. Tardos, Local chromatic number and topological
  properties of graphs, {\em DMTCS Proceedings Series}, {\bf AE} (2005),
  Proceedings of the  European Conference on Combinatorics, Graph Theory and
  Applications, 375--378.  

\bibitem{Stieb} M. Stiebitz, Beitr\"age zur Theorie der f\"arbungskritischen
Graphen, Habilitation, TH Ilmenau, 1985.

\bibitem{Tar} C. Tardif, Fractional chromatic numbers of cones over graphs,
{\em J. Graph Theory}, {\bf 38} (2001), 87--94.

\bibitem{Walk} J.\ W. Walker, From graphs to ortholattices and equivariant
  maps, {\em J. Combin. Theory Ser. B}, {\bf 35} (1983), 171--192.

\bibitem{Y} D.\,A. Youngs, $4$-chromatic projective graphs, {\em J. Graph
Theory}, {\bf 21} (1996), 219--227. 

\bibitem{Ziv} R.\,T. \v{Z}ivaljevi\'c, $WI$-posets, graph complexes and
$\2$-equivalences, {\em J. Combin. Theory Ser. A}, {\bf 111} (2005),  no. 2,
204--223, arXiv:math.CO/0405419.

\bibitem{Zivketlap} R.\,T. \v{Z}ivaljevi\'c, personal communication.

\end{thebibliography}
\end{document}